\documentclass[journal]{IEEEtran}

\usepackage{color}
\usepackage[square,numbers,sort&compress]{natbib}
\usepackage{amssymb}
\usepackage{amsmath}
\usepackage{amsthm}
\usepackage{balance}
\usepackage{epsfig}
\usepackage{eurosym}
\usepackage{graphicx}
\usepackage{graphics}
\usepackage{float}
\usepackage{subcaption}
\usepackage{xspace}

\newtheorem{assumption}{Assumption}

\newtheorem{remark}{Remark}
\newtheorem{definition}{Definition}
\newtheorem{corollary}{Corollary}
\newtheorem{theorem}{Theorem}
\newtheorem{proposition}{Proposition}

\DeclareMathOperator*{\argmin}{arg\,min}
\def\etal{\emph{et al.\, }}

\begin{document}

\title{Distributionally Robust Trading Strategies for \\ Renewable Energy Producers}
%
%
%

\author{Pierre Pinson,~\IEEEmembership{Fellow,~IEEE}
\thanks{P. Pinson has primary affiliation with the Dyson School of Design Engineering, Imperial College London, London, SW7 2AZ, United Kingdom (e-mail: p.pinson@imperial.ac.uk), as well as a secondary affiliation with the Department of Technology, Management and Economics at the Technical University of Denmark, 2800 Kongens Lyngby, Denmark.\\ The research leading to this work is being carried out as a part of the Smart4RES project (European Union’s Horizon 2020, No. 864337). The sole responsibility of this publication lies with the author.
The European Union is not responsible for any use that may be made of the information contained
therein. }
\thanks{Manuscript received xxx; revised xxx.}}

%
%

\markboth{IEEE Transactions, manuscript under review}%
{Pinson: Distributionally Robust Trading for Renewable Energy Producers}
%



\maketitle

\begin{abstract}
Renewable energy generation is offered through electricity markets, quite some time in advance. This then leads to a problem of decision-making under uncertainty, which may be seen as a newsvendor problem. Contrarily to the conventional case for which underage and overage penalties are known, such penalties in the case of electricity markets are unknown, and difficult to estimate. In addition, one is actually only penalized for either overage or underage, not both. Consequently, we look at a slightly different form of a newsvendor problem, for a price-taker participant offering in electricity markets, which we refer to as Bernoulli newsvendor problem. After showing that its solution is consistent with that for the classical newsvendor problem, we then introduce distributionally robust versions, with ambiguity possibly about both the probabilistic forecasts for power generation and the chance of success of the Bernoulli variable. Both versions of the distributionally robust Bernoulli newsvendor problem admit closed-form solutions. We finally use simulation studies, as well as a real-world case-study application, to illustrate the workings and benefits from the approach.
\end{abstract}


%
\IEEEpeerreviewmaketitle

\section{Introduction}
\label{sec:DRintro}

The decarbonization of energy systems, combined with the liberalization of energy markets, makes renewable energy generation increasingly present in electricity markets. In some countries and areas of the world, renewable energy is already reaching very significant shares of the supply. Especially, wind and solar energy are seen as renewable energy sources that could become the major forms of energy generation in many parts of the world. However, owing to their variability in power output, non-dispatchable nature, and limited predictability, wind and solar energy are also bringing challenges in electricity markets. An extensive overview of those aspects are covered in~\cite{Morales2014} among others.

When concentrating on renewables in electricity markets, many have looked at approaches to rethink electricity markets, e.g. in terms of design and pricing, to better accommodate renewable energy generation and its specifics, hence taking a system's view. However, many have also investigated the participation of renewables in existing electricity markets, hence following the agent's perspective. Already when analysing the characteristics of the regulation market in the Nord Pool in 1999, it was established that the asymmetry of regulation penalties (i.e., the spread between forward and balancing prices) could incentivize strategic behaviour from renewable energy producers \cite{Skytte1999}. A number of works followed that looked at various ways to exploit this asymmetry in regulation penalties within various markets, see e.g., \cite{Bathurst2002, Matevosyan2006}. This is while others took a forecasting angle to that problem, by aiming to show that the optimal value of renewable energy in electricity markets would be obtained by using ensemble forecasts \cite{Roulston2003} and probabilistic forecasts of renewable energy generation \cite{Bremnes2004}. Importantly, this problem of market participation for renewable energy generation was recognized as a newsvendor problem \cite{Pinson2007}, and also placed in a general stochastic programming framework \cite{Morales2009}. Additional interesting analytical results were offered for the case of a price-maker setup \cite{Dent2011} and for the specific case of the one-price imbalance settlement (as in the UK) \cite{Browell2018}, among others. The literature on renewables in electricity markets is now expanding rapidly with tens, if not hundreds, of papers being published every year.

Most works focusing on the participation of renewable energy generation in electricity markets rely on the assumption such that those who offer act as price-taker, i.e., that their decisions do not impact market outcomes (both in terms of prices and volumes). \textcolor{black}{We also consider a price-taker setup here, which corresponds to the case of nearly all renewable energy producers in electricity markets. For the more general case of renewable energy producers within a price-maker setup, the reader is referred to the example work of \cite{Dent2011, Zugno2013, Kakhbod2021}.} However, the literature has placed little focus on the fact that, in contrast to the classical newsvendor setup, the actual underage and overage penalties are not known. In most works, it is simply assumed that these may be estimated based on statistics (using historical data), or possibly that these may be predicted with statistical and machines learning approaches. However, it is clear that obtaining high-quality estimates and forecasts is difficult, and that the used underage and overage penalties may substantially deviate from the current conditions. More generally, the probabilistic forecasts used as input are also not perfect, and the true distribution for the uncertain generation may deviate from the predicted one. Consequently, we propose here a distributionally robust version of the newsvendor problem, inspired by renewable energy offering in electricity markets. Considering the general case of decision-making under uncertainty with imperfect forecasts and imperfect estimates of loss functions, Ref.~\cite{vanParys2021} proved that optimal decisions are obtained through distributionally robust optimization. More specifically for the problem of interest here, existing works on distributionally robust newsvendor problems usually focus on the uncertain parameter (being demand or production), not on the underage and overage penalties. Recent relevant examples include the works in~\cite{Fu2021}, \cite{Lee2021} and \cite{Rahimian2019} (among others), who all provide closed-form solutions to the distributionally robust newsvendor problem, but for which ambiguity sets are for the uncertain production (or demand) only. \textcolor{black}{In addition, Ref.~\cite{Rahimian2019} clearly describes the connection between distributionally robust optimization and optimization of coherent risk measures (e.g., conditional Value at Risk -- cVaR).} Recent literature reviews \cite{Rahimian2019b, Lin2022} on distributionally robust optimization give a good introduction to the topic, while pointing at relevant methodological and application-related challenges.

With that objective in mind, we first describe the Bernoulli newsvendor problem, i.e., for which the asymmetry between overage and underage penalties follows a Bernoulli distribution. It comprises a straightforward generalization of the classical newsvendor problem. It is also consistent with the problem of renewable energy producers offering in electricity markets with a two-price imbalance settlement. We look here at distributionally robust Bernoulli newsvendor problems where ambiguity may be about {\it (i)} the probabilistic forecasts for the uncertain parameter, or {\it (ii)} the chance of success of the Bernoulli variable. For the first case, we recall and use the approach introduced in~\cite{Fu2021}. For the second case, we introduce an original closed-form solution, with related proof. 

The paper is structured as following: Section~\ref{sec:preelim} introduces the Bernoulli newsvendor problem and its connection to offering for renewable energy producers in electricity markets. Section~\ref{sec:DRnewsvendor} concentrates on the distributionally robust Bernoulli newsvendor problem, where ambiguity sets are defined for the probabilistic forecasts for the uncertain parameter and for the chance of success of the Bernoulli variable. Section~\ref{sec:simul} provides a set of simulation results for stylized versions of the problem, in order to illustrate the workings of the approach. Eventually, Section~\ref{sec:appl} gathers application results for a real-world case study using wind power generation and electricity market data from France. Finally, Section~\ref{sec:concl} gathers conclusions and perspectives for future works.

\section{Preliminaries: Newsvendor Problems and Electricity Markets}
\label{sec:preelim}

\subsection{Electricity Market Framework}

Consider an electricity market with both forward and balancing stages. The forward stage commonly is day-ahead, and the balancing stage close to real-time. A typical example would be that of the Nord Pool, where all participants place offers before noon for hourly program time units ranging from midnight to midnight the following day. Decisions have to be made at a given time $t$ (before market gate closure) for a set of lead times $\{t+k\}_k$ in the future, where those lead times correspond to the market time units. To simplify our framework and lighten notations, time indices are discarded. This implicitly relies on the assumption that decisions on particular days and given program time units can be made independently of each other. This assumption can be deemed reasonable since there are no inter-temporal constraints involved.

At the forward stage, the market takes the form of a discrete double-sided auction, where both suppliers and consumers place their anonymous offers, which are matched at once following a social-welfare maximization principle. It is said to be discrete since it covers discrete program time units, most often with hourly resolution today in Europe. Out of the market clearing come the production and consumption schedules for all market participants, as well as the equilibrium price $\pi_s$, for each and every program time unit. Under uniform pricing all scheduled energy consumption is bought at $\pi_s$ and all scheduled energy supply is paid at $\pi_s$. 

Following the balancing stage, all deviations from schedule are settled based on the balancing price $\pi_b$, and possibly the state of the system as a whole. Indeed, if considering a one-price imbalance settlement, it is directly $\pi_b$ that is used for settling all imbalances. And, if considering a two-price imbalance settlement instead, only those whose imbalances contribute to the overall system imbalance are subjected to $\pi_b$. This while those who actually help restoring system balance (since their own imbalance is opposite to the system imbalance) are subjected to $\pi_s$. \textcolor{black}{We denote the \emph{overall} system imbalance (also referred to as system length, hence considering the position of all agents in the market) as $s_L$: it is positive if supply is greater than demand overall, and negative if demand is greater than supply overall.}

\textcolor{black}{We place emphasis on two-price imbalance settlement, since this case is more general than the one-price imbalance settlement case. Actually, the solution for the latter case is directly connected to that for the former one, as shown in  Ref.~\cite{Browell2018}: offering strategies are then binary (all or nothing), depending on whether imbalances may yield a reward or a penalty, in expectation. And, the threshold to decide on the binary outcome is given by the same quantity that determines the optimal quantile for the two-price imbalance settlement case. Similar extensions for distributionally robust strategies adapted to the one-price imbalance settlement case could be derived in the future.}

\textcolor{black}{There has been continuous discussion about the merits and caveats of both imbalance settlement approaches over the years. The two-price imbalance settlement is to be preferred if aiming to incentivize truthfulness at the day-ahead stage, since maximum revenue is obtained from null imbalance. In contrast, the one-price imbalance settlement may yield strategic behaviour since it is possible to get a higher revenue when generating imbalances (in the right direction). However, it is then argued that this may turn into a benefit for the system, since market participants who help the system get back to balance will be rewarded. As of today, there has been a shift towards one-price imbalance settlement in most electricity markets. This may change again in the future, depending on how market participants behave in the market, as well as on consequences in terms of price dynamics and imbalances to be managed by system operators.}

\subsection{Renewable Energy Offering as a Bernoulli Newsvendor Problem}

For renewable energy producers, there is uncertainty in future energy generation, since it is not dispatchable and has limited predictability. Therefore, at the time $t$ of placing an offer at the forward stage for a given program time unit $t+k$, energy generation is seen as a random variable, for which a probabilistic forecast is made available. Since the temporal setup is fairly straightforward, we then do not use time-related subscripts in the following, in order to lighten notations.

The random variable $\omega$, for renewable energy generation at that program time unit, takes value in $[0,1]$. It is within that range since it is scaled by the nominal capacity of the power generation asset at hand. The cumulative distribution function (c.d.f.) is denoted $F_\omega$, and the probability density function (p.d.f.) $f_\omega$. The renewable energy producer is to place an offer $y$, as a result of decision-making under uncertainty. Eventually, the energy generation $\omega
^*$ is observed, and we obtain a realization of the revenue of the renewable energy producer. This revenue hence is a random variable, formally defined in the following. 

\begin{definition}[revenue of a renewable energy producer]
The revenue of the renewable energy producer, for an offer $y$ and any value of the uncertain energy generation $\omega$, is given by
\begin{subequations}
\begin{equation}
\label{eq:revenue}
 \mathcal{R}(y,\omega,\pi_s,\pi_b) = \pi_s y + \tilde{\pi}_b (\omega - y) \ ,
\end{equation}
where, considering a two-price imbalance settlement, one has
\begin{equation}
    \tilde{\pi}_b = \left\{ \begin{array}{l} 
    \pi_b \ , \quad   (\omega - y) s_L >0 \\
    \pi_s \ , \quad  (\omega - y) s_L \leq 0 \, .
    \end{array} \right.
\end{equation}
\end{subequations}
\end{definition}

\textcolor{black}{With the above formulation, we overlook the case where there is no energy balancing needed at the second stage, since in that case, $\tilde{\pi}_b = \pi_s$ whatever happens. Hence, renewable energy producers are not penalized for any potential deviation from their contract $y$. Consequently, it means that renewable energy producers will get their optimal revenue whatever happens (even if they bid any random value), as if they could have used an oracle offering strategy allowing to offer a quantity that is exactly what they are to produce.}

In practice, based on the workings of electricity markets, we necessarily have that
\begin{subequations}
\begin{align}
    \tilde{\pi}_b \leq \pi_s, & \quad s_L>0 \, , \\
    \tilde{\pi}_b \geq \pi_s, & \quad s_L<0 \, ,
\end{align}
\end{subequations}
which makes that we can define the following overage and underage penalties, as a basis to describe a newsvendor problem.

\begin{definition}[overage and underage penalties]
Based on forward and balancing prices $\pi_s$ and $\pi_b$, as well as the overall system state $s_L$, overage $\pi_o$ and underage $\pi_u$ penalties can be defined
\begin{subequations}
\begin{align}
    \label{eq:overage}
    \pi_o & \textcolor{black}{= (\pi_s - \tilde{\pi}_b) = } \, (\pi_s - \pi_b) \, \mathbf{1}_{\{ s_L \geq 0 \}} \ , \\
    \pi_u & \textcolor{black}{= (\tilde{\pi}_b - \pi_s) = } \, (\pi_b - \pi_s) \, \mathbf{1}_{\{ s_L < 0 \}} \ ,
    \label{eq:underage}
\end{align}
\end{subequations}
where $\mathbf{1}_{\{ . \}}$ is the indicator function, hence returning 1 if the statement ${\{ . \}}$ is true, and 0 otherwise. 
\end{definition}

\textcolor{black}{In most of the literature covering offering strategies for renewable energy producers, as well as newsvendor problems, underage and overage penalties are assumed to be known (and hence, deterministic). This is not the case in the following: the system length $s_L$ is a random variable and the price deviation $\pi_s - \pi_b$ is also a random variable. Consequently, the overage $\pi_o$ and underage $\pi_u$ penalties are random variables too.}

\textcolor{black}{After a little algebra (for instance covered in~\cite{Pinson2007}), starting from the revenue function~\eqref{eq:revenue}, we want to obtain the corresponding opportunity cost function to be minimized. We can first rewrite~\eqref{eq:revenue} as if a loss function to be minimized,
\begin{equation}
\label{eq:lossderiv1}
 - \pi_s y - \tilde{\pi}_b (\omega - y) \ .
\end{equation}
The second term can be written as the sum of 2 sub-cases for when \emph{(i)} generation is greater than the decision $y$, and when \emph{(ii)} generation is less than decision $y$,
\begin{equation}
\label{eq:lossderiv2}
 - \pi_s y - \tilde{\pi}_b (\omega - y)_+ +  \tilde{\pi}_b (y-\omega)_+ \ ,
\end{equation}
where $(.)_+$ is for the positive part, that is, $u_+ = u \,  \mathbf{1}_{u\geq0}$. We then introduce 2 terms that sum to 0, i.e.,
\begin{equation}
\label{eq:lossderiv3}
 - \pi_s y \, \underbrace{ + \,  \pi_s \omega - \pi_s \omega}_{=0}  - \tilde{\pi}_b (\omega - y)_+ +  \tilde{\pi}_b (y-\omega)_+\ .
\end{equation}
By gathering terms, we obtain
\begin{equation}
\label{eq:lossderiv5}
 - \pi_s \omega - \pi_s (y-\omega) - \tilde{\pi}_b (\omega - y)_+ +  \tilde{\pi}_b (y-\omega)_+\ .
\end{equation}
The $\pi_s (y-\omega)$ term can eventually be distributed to the 2 sub-cases, as
\begin{equation}
\label{eq:lossderiv5}
- \pi_s \omega + (\pi_s - \tilde{\pi}_b) (\omega - y)_+ + (\tilde{\pi}_b-\pi_s) (y-\omega)_+\ .
\end{equation}
In the above, we recognize the definition of overage and underage penalties as in~\eqref{eq:overage} and~\eqref{eq:underage}. In addition, the first term $- \pi_s \omega$ is not a function of the decision $y$. Hence, in terms of a loss function to be minimized, the optimal decision will be the same if the term is removed. This comes in contrast with the formulation in \cite{Rahimian2019}, though aligned with that in~\cite{Fu2021}. Similarly, that loss function can be divided by a term that is independent of the decision $y$ and therefore not affect its use for optimization purposes. This eventually yields the final loss function:
\begin{equation} \label{eq:scaledopcost}
    \mathcal{L}(y,\omega,\pi_o,\pi_u) = \frac{\pi_o}{\pi_o+\pi_u} (\omega-y)_+ + \frac{\pi_u}{\pi_o+\pi_u}  (y-\omega)_+ \, . 
\end{equation}}

\textcolor{black}{Compared to the conventional loss function for newsvendor problems, for which underage and overage costs are known, they are here 2 inter-related random variables (since their sum is equal to 1). If defining the first random variable as $s = \frac{\pi_o}{\pi_o+\pi_u}$, the second one is $1-s$. $s$ necessarily is a Bernoulli random variable since it has only 2 potential discrete outcomes, i.e.
\begin{align}
    (i) & \quad \pi_o \neq 0, \, \, \pi_u = 0 \, \,  \Rightarrow\, \, s=\frac{\pi_o}{\pi_o+\pi_u}=1, \, \, 1-s=0 \nonumber \\
    (ii) & \quad  \pi_o = 0, \, \, \pi_u \neq 0 \, \, \Rightarrow \, \, s=\frac{\pi_o}{\pi_o+\pi_u}=0, \, \, 1-s=1 \nonumber 
\end{align}
Eventually, this means that we can reduce all uncertainties in \emph{both} system length and underage/overage penalties to a single Bernoulli random variable with chance of success $\tau$, $s \sim \text{Bern}(\tau)$. And, the chance of success $\tau$ may be seen as the expected asymmetry between overage and underage costs, accounting for system length, i.e., $\tau = \mathbb{E} \left[\pi_o/(\pi_o+\pi_u) \right]$.}

In that framework, the price-taker assumption involves simplifications in terms of the dependencies between $\omega$, $s$ and the decision variables $y$.

\begin{assumption}[price-taker assumption]
For a Bernoulli newsvendor problem, the price-taker assumption  translates to
\begin{itemize}
\item[(A1)] $\omega$ and $s$ are independent random variables,
\item[(A2)] the distributions of $\omega$ and $s$ are independent of the decision $y$.
\end{itemize}
\end{assumption}

As for any other approach to the participation of renewable energy producers in electricity markets, one may wonder whether participation strategies relying on a price-taker assumption could not actually result in affecting clearing outcomes if those are employed by a non-negligible population of renewable energy producers. This may then lead to more complex strategies accounting for individual impact on market outcomes, e.g., \cite{Zugno2013}, or accounting for population impact on market outcomes, e.g., \cite{Kakhbod2021}. In all cases, this is not specific to employing distributionally robust optimization and a research question by itself. \textcolor{black}{In addition, Assumption A1 implies that the wind power generation of a single wind farm (or portfolio) offering in the market does not impact the balance of the system. If wind energy in a given system was the main driver of overall system balance, this would translate to saying that $\omega$ is weakly correlated to overall wind power generation in the system. However, in practice, system balance is induced by more potential causes than wind power generation alone, and the fact that $\omega$ and $s$ are not correlated can be deemed reasonable.}

Consequently, we formally define the Bernoulli newsvendor problem, to be used as a basis for the offering of renewable energy producers in electricity markets.

\begin{definition}[Bernoulli newsvendor problem]
Based on a Bernoulli random variable $s$ (with chance of success $\tau$), and the uncertain parameter $\omega$ (with c.d.f.\ $F_\omega$), the decision $y^*$ minimizing the expected opportunity cost function is
\begin{subequations}
\begin{equation} \label{eq:Bernnews}
    y^* = \argmin_{y} \, \mathbb{E}_{\omega,s} \left[ \mathcal{L}(y,\omega,s) \right] \, ,
\end{equation}
where the opportunity cost is defined as
\begin{equation} \label{eq:scaledopcost2}
    \mathcal{L}(y,\omega,s) = s (\omega-y)_+ + (1-s) (y-\omega)_+ \, ,
\end{equation}
\end{subequations}
\textcolor{black}{and where $(.)_+$ is for the positive part.}
\end{definition}

The Bernoulli newsvendor problem has a solution which is readily connected to that for the classical newsvendor problem, except that, instead of relying on the ratio of overage and sum of penalties, it involves the chance of success of the Bernoulli variable.

\begin{proposition}
\label{prop:Bnewsvendor}
Consider $F_\omega$ the c.d.f.\ for the uncertain parameter $\omega$ and $\tau$ the chance of success for $s$. The optimal decision $y^*$ for the Bernoulli newsvendor problem~\eqref{eq:Bernnews} is
\begin{equation}
    y^* = F_\omega^{-1} (\tau) \, .
\end{equation}
\end{proposition}

The proof of Proposition~\ref{prop:Bnewsvendor} is given in Appendix~\ref{app:proofnewsvendor}. The above result means that, as for the classical newsvendor problem where the predictive distribution and penalties both are considered as known, the optimal offer is a specific quantile of the predictive c.d.f.\ for the uncertain parameter $\omega$. The nominal level of that quantile is $\tau$, the chance of success of the Bernoulli variable. In the following, since the Bernoulli variable $s$ is uniquely defined by its chance of success $s$, we use the notation $\mathcal{L}(y,\omega,s)$ and $\mathcal{L}(y, \omega, \tau)$ interchangeably. 

In practice, a model is proposed to predict the chance of success for the Bernoulli variable, to be used as a basis to obtain the optimal quantile. The easiest strategy would be to compute average values of overage and underage penalties over a period with historical data, as well as frequencies for the system being long and short, to then deduce an estimate (also seen as a forecast) $\hat{\tau}$ for the chance of success. Similarly, the cumulative distribution function $F_\omega$ is not readily available and needs to be predicted. We write $\hat{F}_\omega$ that predictive c.d.f., most often provided by expert forecasters based on weather forecasts and statistical or machine-learning approaches.

\section{Distributionally Robust Newsvendor Problem}
\label{sec:DRnewsvendor}

While in the classical newsvendor setup, it is assumed that the predictive distribution for the uncertain parameter $\omega$, as well as the overage and underage penalties, are known, it is not the case here. This also means that the distribution of $s$ (fully characterized by the chance of success $\tau$) is not known. As input to decision-making, one has a predictive c.d.f. $\hat{F}_\omega$ for $\omega$ and a predictive c.d.f. $\hat{F}_s$ for the Bernoulli variable $s$, which is readily given by the forecast $\hat{\tau}$ for the chance of success. We expect the true distribution for $\omega$ and for $s$ to deviate from $\hat{F}_\omega$ and $\hat{F}_s$, respectively. Distributionally robust optimization allows us to accommodate that ambiguity. In the following, we first recall the solution proposed by Fu \etal~\cite{Fu2021}, when the ambiguity is about $\hat{F}_\omega$. We subsequently describe our original solution for the case where the ambiguity is about $\hat{F}_s$. 

\textcolor{black}{
\begin{remark}
While in the expected utility maximization case (as in Proposition~\ref{prop:Bnewsvendor}), it is straightforward to assume that the term $\pi_s \, \omega$ can be overlooked in order to obtain the final loss function~\eqref{eq:scaledopcost}, this should be reconsidered in the distributionally robust case. This is since, in principle, the term $\pi_s \, \omega$ could affect the worst-case distribution to be considered. However, if considering ambiguity about $\hat{F}_s$, this term can again be overlooked since the worst-case distribution does not involve $\omega$. And, if focusing on ambiguity about $\hat{F}_\omega$, we will explain when stating our main result in the following (Theorem~\ref{theo:DRBnewsvendor1}) why the worst-case distribution is the same whether considering the term $\pi_s \, \omega$ or not.
\end{remark}
}

\subsection{Ambiguity about $\hat{F}_\omega$}

A distributionally robust optimization view of the Bernoulli newsvendor problem takes the form of a min-max optimization problem. If the ambiguity is on the distribution $\hat{F}_\omega$ of the uncertain parameter $\omega$, the worst case is in terms of potential distributions within a certain distance from $\hat{F}_\omega$. And, since forecasts for renewable energy generation most often take a non-parametric form, it is more convenient to consider ambiguity sets defined in terms of distance instead of moments. We write $\mathcal{B}_{\hat{F}_\omega} (\rho)$ that ambiguity set, with radius $\rho$. Let us then formally define that problem in the following.

\begin{definition}[distributionally robust Bernoulli newsvendor problem -- ambiguity about $\hat{F}_\omega$]
Consider a Bernoulli random variable $s$ with estimated chance of success $\hat{\tau}$, the uncertain production $\omega$ with predictive c.d.f.\ $\hat{F}_\omega$, and an ambiguity set $\mathcal{B}_{\hat{F}_\omega} (\rho)$ with radius $\rho$. The distributionally robust Bernoulli newsvendor problem, with ambiguity about $\hat{F}_\omega$, is that for which the decision $y^*$ is given by 
\begin{equation} \label{eq:DRnews1}
    y^* = \argmin_{y} \, \sup_{F_\omega \in \mathcal{B}_{\hat{F}_\omega} (\rho)} \, \mathbb{E}_{\omega,s} \left[ \mathcal{L}(y,\omega,\tau) \right]  \, .
\end{equation}
\end{definition}

For such a problem, where the c.d.f.\ of $\omega$ has bounded support ($\text{supp}(F_\omega) = [\underline{\omega}, \overline{\omega}]$), Fu \etal \cite{Fu2021} showed that it is sufficient to consider two distributions $\underline{F}_\omega$ and $\overline{F}_\omega$ defined based on a $\phi$-divergence and the associated so-called Jager-Wellner discrepancy measure, to define the worst case for the inner problem. This concept for which these 2 distributions define the ambiguity set for $\hat{F}_\omega$ is referred to as a first-order stochastic dominance ambiguity set (FSD-ambiguity set). Eventually, it simply means that
\begin{equation*}
    \underline{F}_\omega(x) \leq F_\omega (x) \leq \overline{F}_\omega(x), \quad \forall x, \, \forall F_\omega \in \mathcal{B}_{\hat{F}_\omega} (\rho)  \, .
\end{equation*}

In a more general manner, this involves two continuous deformation operators (upper $\overline{\mathcal{O}}_\rho$ and lower $\underline{\mathcal{O}}_\rho$) that can alter the c.d.f.\ from its original form $\hat{F}_\omega$ (for the case of $\rho=0$) to the Heaviside functions $H(.)$ located at the bounds of $\text{supp}(F_\omega)$. Let us generally define such deformation operators here.

\begin{definition}[Upper and lower deformation operators] \label{def:deformationoperators}
Consider a reference c.d.f.\ $F_\omega$. Upper $\overline{\mathcal{O}}_\rho$ and lower $\underline{\mathcal{O}}_\rho$ deformation operators are continuous operators on $F_\omega$ such that
\begin{align}
    \text{(identity): } &  \quad \overline{\mathcal{O}}_0 (F_\omega) \, = \, \underline{\mathcal{O}}_0 (F_\omega) \, = \, F_\omega  \nonumber \\
    \text{(robustness): } & \quad \lim_{\rho\rightarrow 1} \overline{\mathcal{O}}_\rho (F_\omega) \, = \, H(0), \, \, \text{and} \nonumber \\
    & \quad \lim_{\rho\rightarrow 1} \underline{\mathcal{O}}_\rho (F_\omega) \, = \, H(1) \nonumber \\ 
    \text{(monotonicity): } & \quad \overline{\mathcal{O}}_\rho (F_\omega) \, \leq \, \overline{\mathcal{O}}_{\rho'}, \quad \underline{\mathcal{O}}_\rho (F_\omega) \, \geq \, \underline{\mathcal{O}}_{\rho'}, \nonumber \\ & \quad \forall \rho,\rho' \in [0,1], \, \, \rho \leq \rho' \nonumber 
\end{align}
\end{definition}

As an example, we introduce here a double-power deformation operator that fulfill the above definition.

\begin{definition}[double-power deformation operator]
Consider a reference c.d.f.\ $F_\omega$. The upper $\overline{\mathcal{O}}_\rho$ and lower $\underline{\mathcal{O}}_\rho$ double-power deformation operators are defined as
\begin{subequations}
\begin{align}
    \overline{\mathcal{O}}_\rho ( F_\omega ) & \, = \, \left( 1- ( 1 - F_\omega{}^{\frac{1}{1-\rho}} ) \right)^{1-\rho} \,, \\ 
    \underline{\mathcal{O}}_\rho ( F_\omega )& \, = \,  1 - (1 - F_\omega{}^{\frac{1}{1-\rho}})^{1-\rho} \, ,
\end{align}
\end{subequations}
with $\rho$ the deformation parameter.
\end{definition}

The double-power deformation operator has the additional property of being symmetric, in the sense that the deformation from the original c.d.f. $F_\omega$ to the two limiting distributions $\underline{F}_\omega$ and $\overline{F}_\omega$ is the same in both directions, i.e., $\overline{F}_\omega(x)-F_\omega(x) = F_\omega(x)-\underline{F}_\omega(x), \, \forall x \in [0,1]$ and $\underline{F}^{-1}_\omega(p)-F^{-1}_\omega(p) = F^{-1}_\omega(p)-\overline{F}^{-1}_\omega(p), \, \forall p \in [0,1]$. Other more general (but not necessarily symmetric) deformation operators could be defined, e.g., exponential-Pareto ones, as used for the example of the estimation of Lorenz curves \cite{Sitthiyot2021}.

In the above, the deformation parameter $\rho$ can readily be used as the ball radius for the definition of the ambiguity set, since it defines the two bounding distributions $\underline{F}_\omega$ and $\overline{F}_\omega$ for the ball $\mathcal{B}_{\hat{F}_\omega}$.

\begin{definition}[FSD-ambiguity set for $\hat{F}_\omega$]
Given a ball radius $\rho$, the FSD-ambiguity set for $\hat{F}_\omega$ is defined as $[\underline{F}_\omega, \, \overline{F}_\omega]$ where
\begin{subequations}
\begin{align}
    \underline{F}_\omega & \, = \, \underline{\mathcal{O}}_\rho ( \hat{F}_\omega ) \, ,\\
    \overline{F}_\omega  & \, = \, \overline{\mathcal{O}}_\rho ( \hat{F}_\omega ) \, .
\end{align}
\end{subequations}
\end{definition}

\begin{figure}[!ht]
\centering
    \includegraphics[width=.93\columnwidth]{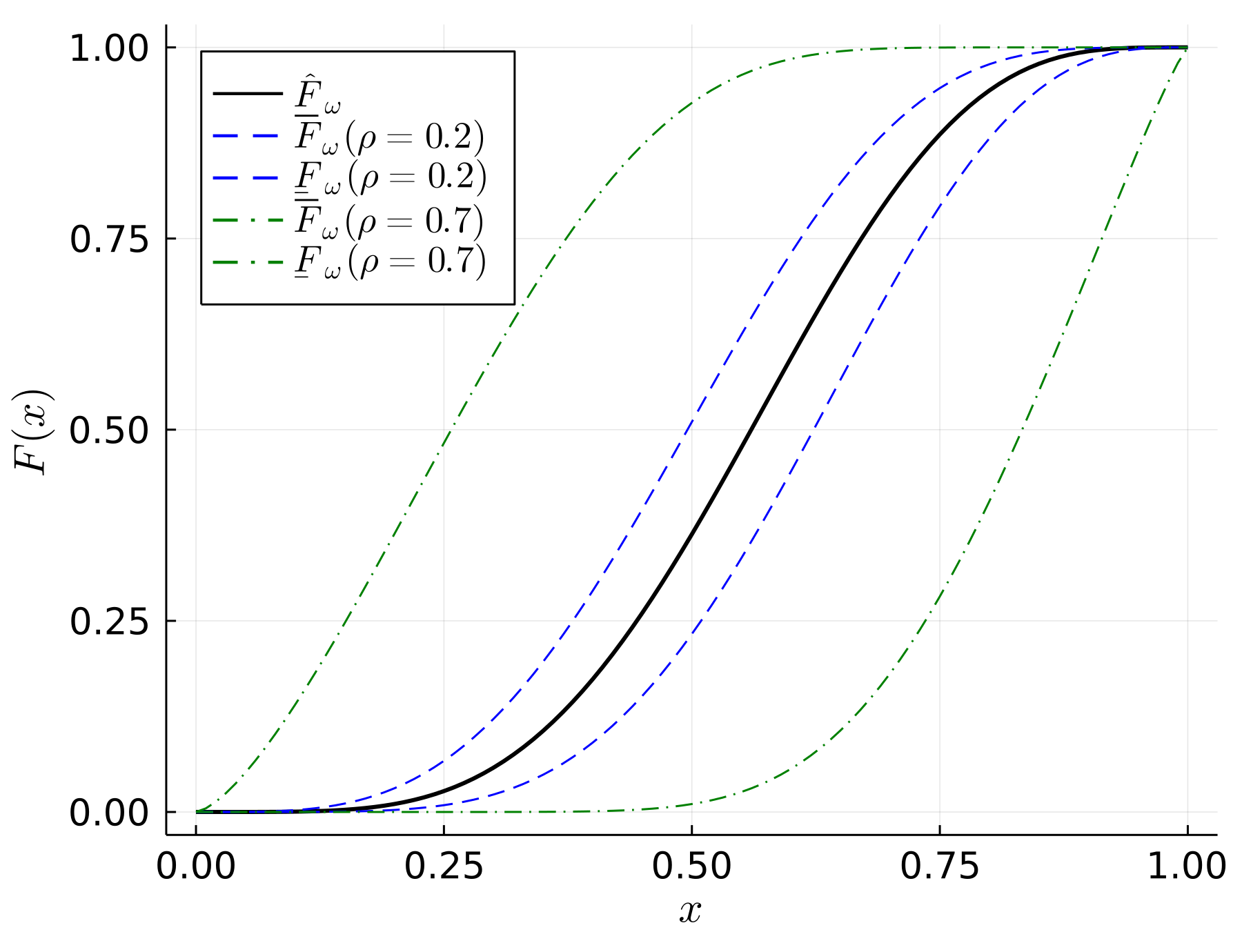}
    \caption{Ambiguity sets for $F_\omega$ obtained using the double-power deformation operator (cases of $\rho=0.2$ and $\rho=0.7$).}
    \label{fig:exdoublepower}
\end{figure}

Figure~\ref{fig:exdoublepower} provides illustrative examples of ambiguity sets (and the related ball $\mathcal{B}_{\hat{F}_\omega}$) obtained by using the double-power deformation operator, with various values for $\rho$. By increasing the ball radius $\rho$, the limiting distributions get further away from the reference one. Importantly, in contrast to the proposals of Fu \etal~\cite{Fu2021}, these do not need to be trimmed at 0 or 1, since whatever $\rho$ , $\overline{\mathcal{O}}_\rho (F_\omega)$ and $\underline{\mathcal{O}}_\rho (F_\omega)$ take values within the unit interval $[0,1]$.

The following theorem gives the solution of the distributionally robust Bernoulli newsvendor problem with ambiguity about $\hat{F}_\omega$. It is adapted after the work of Fu \etal~\cite{Fu2021}. The proof is hence omitted since available in their manuscript.

\begin{theorem}
\label{theo:DRBnewsvendor1}
Consider an FSD-ambiguity set defined by a ball $\mathcal{B}_{\hat{F}_\omega} (\rho)$ with radius $\rho$, yielding the two bounding distributions $\underline{F}_\omega$ and $\overline{F}_\omega$. For a predicted chance of success $\hat{\tau}$, the solution of the distributionally robust Bernoulli newsvendor problem~\eqref{eq:DRnews1} is
\begin{equation} \label{eq:DROsol1}
    y^* \, = \, \hat{\tau} \underline{F}_\omega{}^{-1}(\hat{\tau}) + (1-\hat{\tau})  \overline{F}_\omega{}^{-1}(\hat{\tau}) \, .
\end{equation}
\end{theorem}

\textcolor{black}{
\begin{remark}
The worst-case distribution $F^{\text{ws}}$ for this distributionally robust problem is defined by \begin{equation}
    F^{\text{ws}} (x) = \left\{ \begin{array}{l}
         \overline{F}_\omega(x), \quad x < \overline{F}_\omega{}^{-1}(\hat{\tau}) \\
         \hat{\tau}, \quad \overline{F}_\omega{}^{-1}(\hat{\tau}) < x < \underline{F}_\omega{}^{-1}(\hat{\tau}) \\
         \underline{F}_\omega(x), \quad x > \underline{F}_\omega{}^{-1}(\hat{\tau})
    \end{array} \right. \, .
\end{equation}
Consequently, based on that FSD-ambiguity set, one obtains the same worst case distribution whether considering the term $\pi_s \, \omega$, or not. Indeed, this is because the worst case necessarily is that the probability mass is pushed towards the bounds, hence yielding maximum exposure to overage and underage costs. This is also valid if considering the term $\pi_s \, \omega$ since, when shifting the left part of the distribution to the left (from $\hat{F}_\omega$ to $\overline{F}_\omega$), one jointly minimizes potential revenue $\pi_s \, \omega$ while maximizing potential underage costs. Similarly, when shifting the right part of the distribution to the right (from $\hat{F}_\omega$ to $\underline{F}_\omega$), one jointly minimizes potential revenue $\pi_s \, \omega$ while maximizing potential overage costs.
\end{remark}
}

For sufficiently large values of the radius $\rho$, we obtain the following robust limiting case.

\begin{corollary}
\label{corr:robnewsvendor1}
For sufficiently large values of $\rho$, the radius of $\mathcal{B}_{\hat{F}_\omega} (\rho)$, the solution of the distributionally robust Bernoulli newsvendor problem~\eqref{eq:DRnews1} converges to the robust solution  $y^* \, = \,  \hat{\tau}$.
\end{corollary}

A compact proof of this result is given in Appendix~\ref{app:proofrobustnewsvendor1}.
It is interesting to see that this robust limiting case is equivalent to considering an uninformative predictive distribution $\mathcal{U}[0,1]$ for the uncertain parameter $\omega$. Indeed, if $\hat{F}_\omega$ describes a standard uniform $\mathcal{U}[0,1]$, the optimal quantile for the Bernoulli newsvendor problem is $\hat{F}_\omega^{-1}(\hat{\tau}) = \hat{\tau}$.

\subsection{Ambiguity about $\hat{F}_s$}

While the case of ambiguity about $\hat{F}_\omega$ has already been explored in the literature, there is no result related to ambiguity about $\hat{\tau}$. In that case, the distributionally robust optimization problem also is a min-max problem. However, the worst case is to be defined based on the distribution of $s$.

An appealing feature of that problem is that the distribution of $s$ is straightforward to handle, since it is characterized by single parameter $\tau$ only. Hence, both moment-based and distance-based approach to defining ambiguity sets for $\hat{F}_s$ end up being the same. It translates to considering an interval for $\tau$ to define the ambiguity set, as a ball (an interval, really) around the value $\hat{\tau}$. The ball $\mathcal{B}_{\hat{\tau}}$ around $\hat{\tau}$ is hence given by
\begin{equation}
 \mathcal{B}_{\hat{\tau}} (\varepsilon) = [\underline{\tau}, \, \overline{\tau}], \quad \varepsilon \geq 0 \, .
\end{equation}

There may be different ways to obtain the interval bounds and hence to define the ball $\mathcal{B}_{\hat{\tau}}$, which gives the ambiguity set for $\hat{F}_s$.

\begin{definition}[uniform and level-adjusted ambiguity sets for $\hat{F}_s$]
Given a ball radius $\varepsilon$, a uniform ambiguity set for $\hat{F}_s$ is defined by the ball $\mathcal{B}_{\hat{\tau}}$ with
\begin{subequations}
\begin{align}
    \underline{\tau} &\, = \, \max(\hat{\tau}-\varepsilon, \, 0) \\ 
    \overline{\tau} &\, = \, \min(\hat{\tau}+\varepsilon, \, 1)
\end{align}
\end{subequations}
while a level-adjusted ambiguity set (with parameter $\theta \in [0,1)$) for $\hat{F}_s$ is defined by the ball $\mathcal{B}_{\hat{\tau}}$ with
\begin{subequations}
\begin{align}
    \underline{\tau} &\, = \, \max(\hat{\tau} -  \varepsilon ( 1 \, - \,  4 \, \theta \,  \hat{\tau}(1-\hat{\tau})),\, 0) \\ 
    \overline{\tau} &\, = \, \min(\hat{\tau} \, + \,  \varepsilon ( 1 \, - \,  4 \, \theta \,  \hat{\tau}(1-\hat{\tau})), \, 1)
\end{align}
\end{subequations}
\end{definition}
The interest of the level-adjusted ambiguity set is that is recognizes the fact that the ambiguity may be a function of the chance of success. In practice, it appears reasonable to think that the ambiguity set may be larger closer to the bounds, and smaller for a chance of success $\tau$ close to 0.5.

Now, we can define our distributionally robust Bernoulli newsvendor problem, with ambiguity about $\hat{F}_s$. 

\begin{definition}[distributionally robust Bernoulli newsvendor problem -- ambiguity about $\hat{F}_s$]
Consider a Bernoulli random variable $s$ with estimated chance of success $\hat{\tau}$, the uncertain production $\omega$ with predictive c.d.f.\ $\hat{F}_\omega$, and an ambiguity set for $\hat{F}_s$ defined by the ball $\mathcal{B}_{\hat{\tau}} (\varepsilon)$ with radius $\varepsilon$. The distributionally robust Bernoulli newsvendor problem, with ambiguity about $\hat{F}_s$, is that for which the decision $y^*$ is given by
\begin{equation} \label{eq:DRnews2}
    y^* = \argmin_{y} \, \max_{\tau \in \mathcal{B}_{\hat{\tau}}(\varepsilon)} \, \mathbb{E}_{\omega,s} \left[ \mathcal{L}(y,\omega,\tau) \right]
\end{equation}
\end{definition}

Since having a quite simple setup, it is possible to derive a closed-form solution to this distributionally robust Bernoulli newsvendor problem with ambiguity about $\hat{F}_s$.

\begin{theorem}
\label{theo:DRBnewsvendor2}
Consider an ambiguity set for $\hat{F}_s$ defined by a ball $\mathcal{B}_{\hat{\tau}} (\varepsilon)$ with radius $\varepsilon$, and the predictive c.d.f.\ $\hat{F}_\omega$ for the uncertain parameter $\omega$. The solution of the distributionally robust Bernoulli newsvendor problem~\eqref{eq:DRnews2} is
\begin{align} \label{eq:DROsol2}
    y^* \, = & \, \, \hat{F}_\omega^{-1}(\overline{\tau}) \,   \mathbf{1}_{\{ \hat{F}_\omega^{-1}(\overline{\tau}) < \mathbb{E}[\omega]\}} + \hat{F}_\omega^{-1}(\underline{\tau}) \,  \mathbf{1}_{\{ \hat{F}_\omega^{-1}(\underline{\tau})  >\mathbb{E}[\omega]\}} \nonumber \\
     & + \, \mathbb{E}[\omega] \, \mathbf{1}_{\{ \hat{F}_\omega^{-1}(\overline{\tau}) \geq \mathbb{E}[\omega]\}} \, \mathbf{1}_{\{ \hat{F}_\omega^{-1}(\underline{\tau}) \leq\mathbb{E}[\omega]\}} \, .
\end{align}
\end{theorem}

\begin{figure*}[!ht]
\centering
\begin{subfigure}{0.31\textwidth}
         \includegraphics[width=\columnwidth]{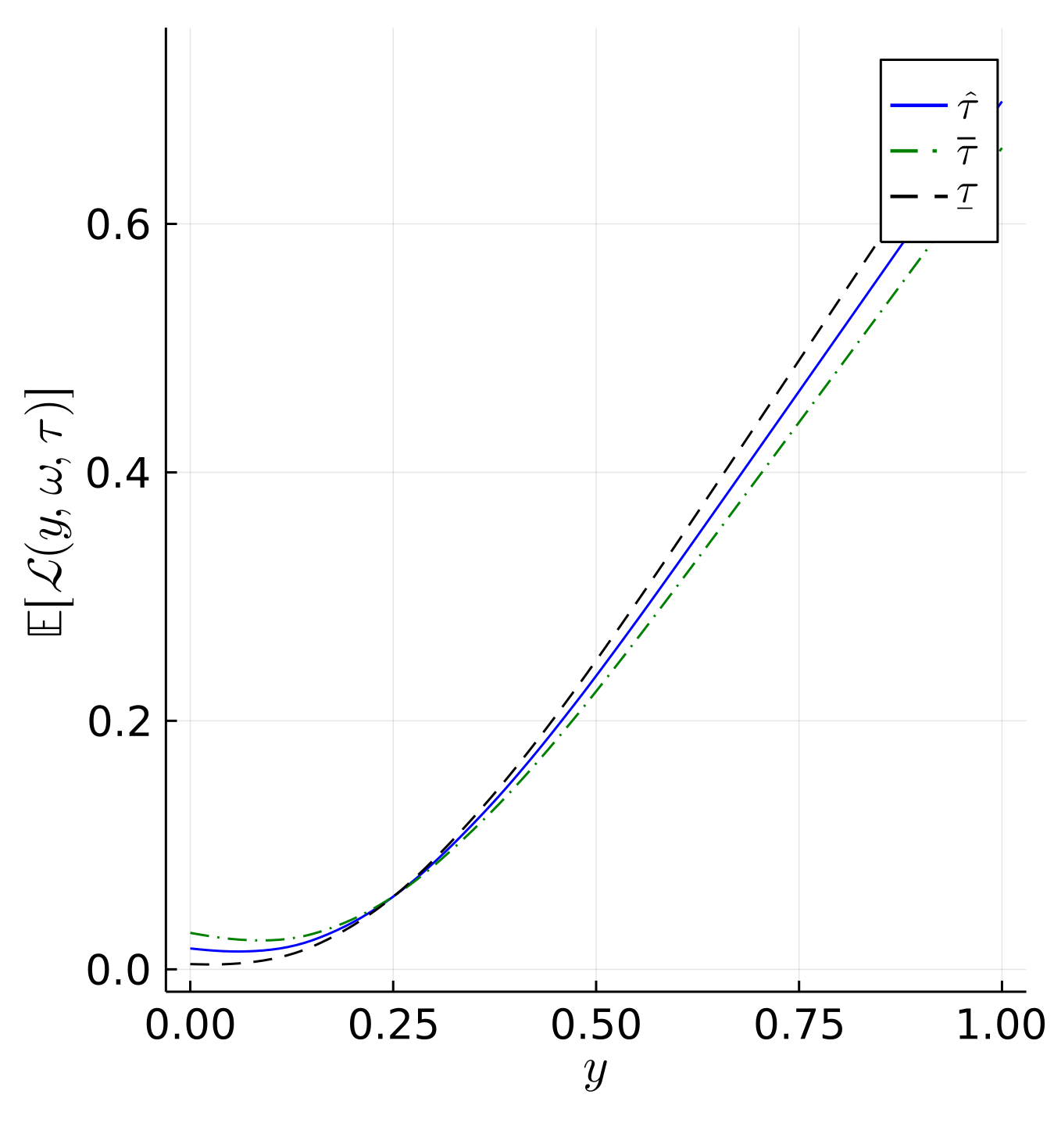}
         \caption{$\tau=0.1$, $\varepsilon=0.05$}
         \label{fig:DRnewsvendor-example1}
\end{subfigure}
\hfill
\begin{subfigure}{0.31\textwidth}
         \includegraphics[width=\columnwidth]{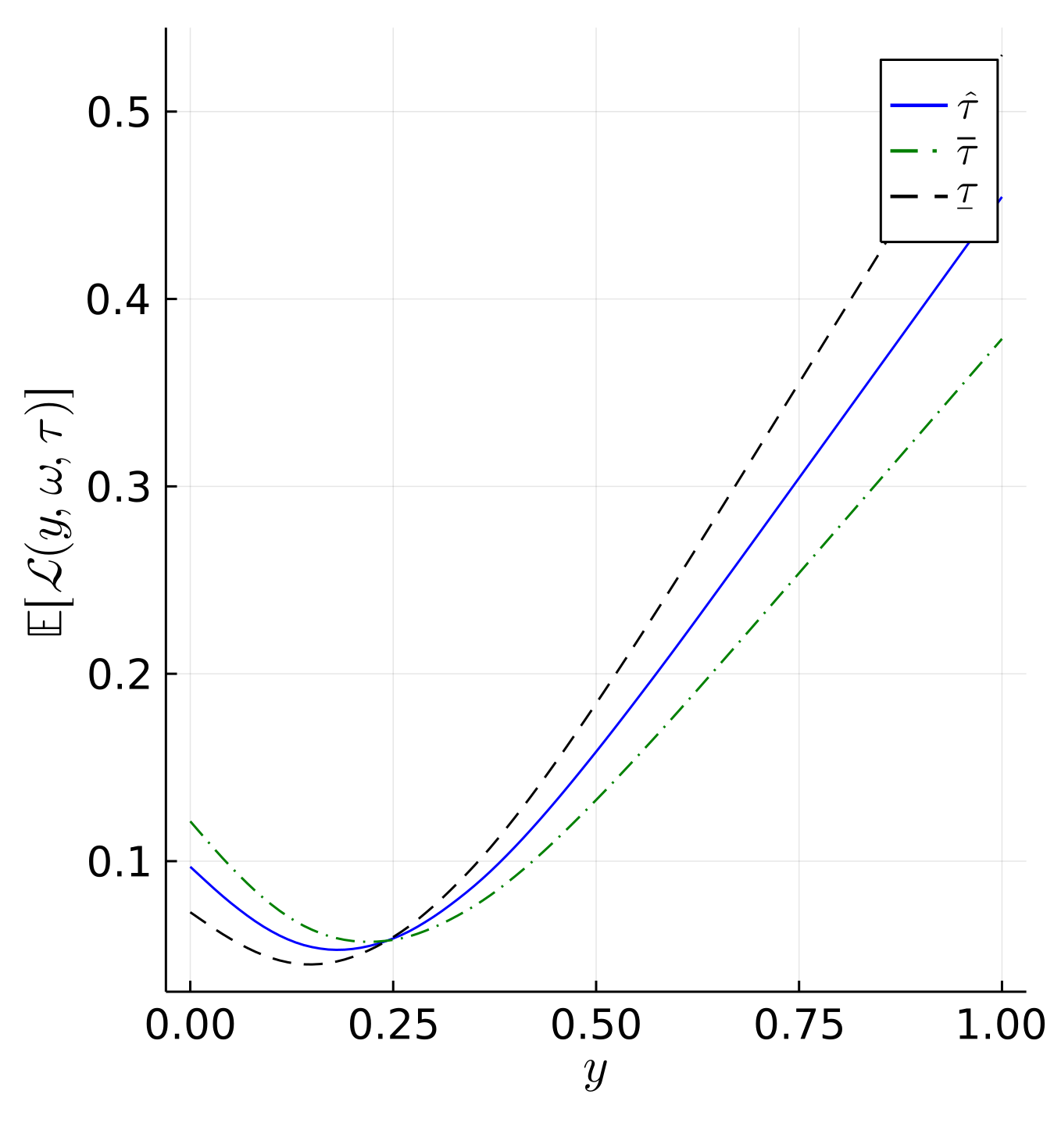}
         \caption{$\tau=0.5$, $\varepsilon=0.1$}
         \label{fig:DRnewsvendor-example2}
\end{subfigure}   
\hfill
\begin{subfigure}{0.31\textwidth}
         \includegraphics[width=\columnwidth]{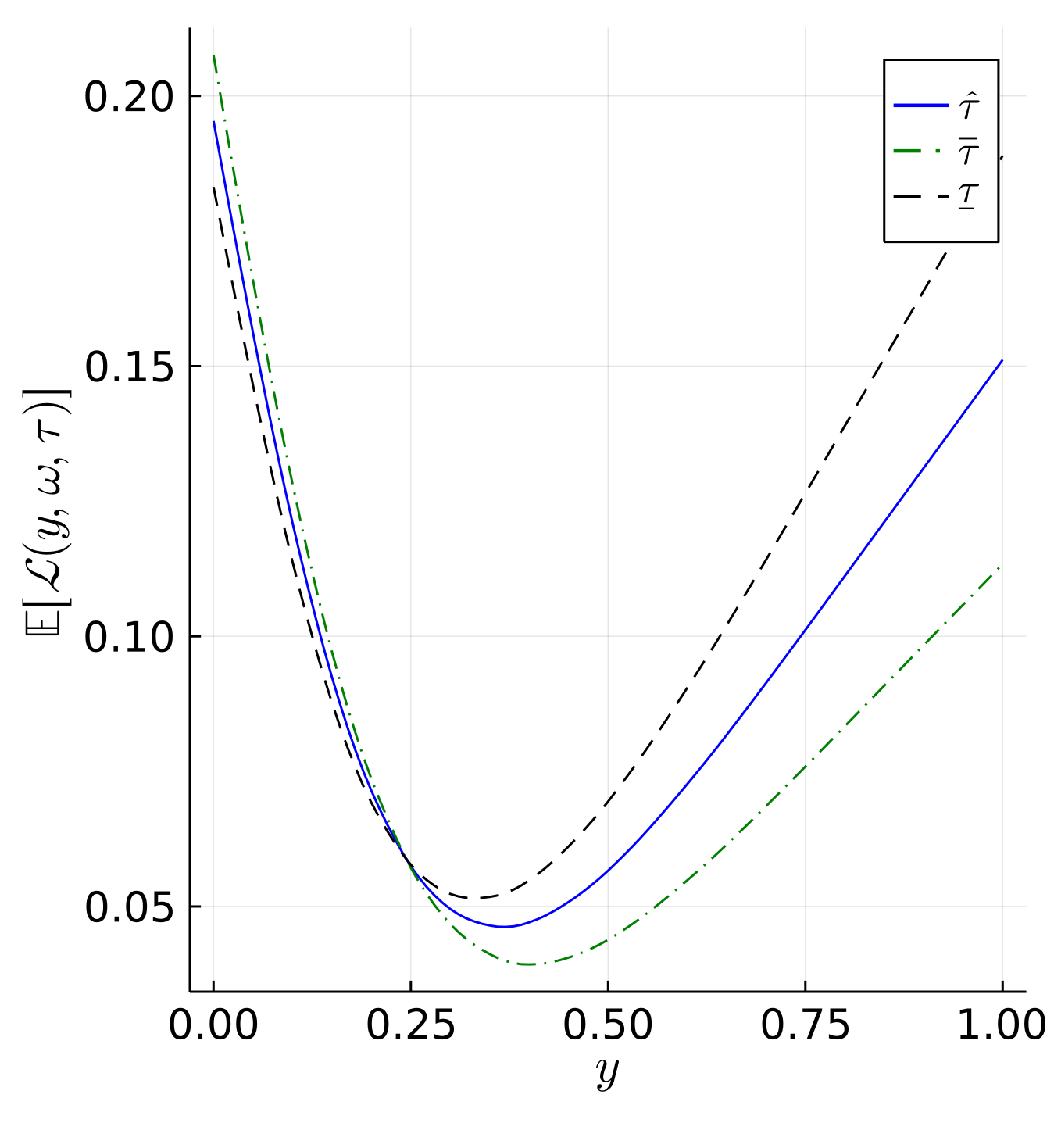}
         \caption{$\tau=0.8$, $\varepsilon=0.05$}
         \label{fig:DRnewsvendor-example3}
\end{subfigure}   
\caption{Expected opportunity cost, as a function the decision $y$, for different values of $\tau$ and $\varepsilon$ (using uniform ambiguity sets). The uncertain parameter $\omega$ follows a \textcolor{black}{Beta(2,6)} distribution \textcolor{black}{(with expected value $\mathbb{E}[\omega]=0.25$)}, while the estimate $\hat{\tau}$ is based on 15 samples. \label{fig:DRnewsvendor-example}}
\end{figure*}

The proof of Theorem~\ref{theo:DRBnewsvendor2} is given in Appendix~\ref{app:proofDRnewsvendor2}. In addition, let us provide in the following an intuition for the result, based on Figure~\ref{fig:DRnewsvendor-example}. \textcolor{black}{We show there the expected opportunity cost as a function of the decision $y$, based on an uncertain parameter $\omega$ following a Beta(2,6) distribution (hence, with potential outcomes in [0,1] and with expected value $\mathbb{E}[\omega] = 0.25$).} \textcolor{black}{The expected loss $\mathbb{E}_\omega \left[ \mathcal{L}(y,\omega,\tau) \right]$ (on the $y$-axis) can be calculated based on~\eqref{eq:scaledopcost} where the expectation $\mathbb{E} \left[\pi_o/(\pi_o+\pi_u) \right]$ is given by $\tau$.} An estimate of $\tau$ is obtained from 15 values sampled from a Bern$(\tau)$ distribution for $s$. \textcolor{black}{Figure 2 gathers 3 illustrative cases, for different chosen values for the chance of success $\tau$ and ball radius $\varepsilon$.} For these 3 cases, we visualize the expected opportunity cost for $\hat{\tau}$, $\underline{\tau}$ and $\overline{\tau}$ for a uniform ambiguity set for $\hat{F}_s$. The first thing to observe is that all curves cross for $y=\mathbb{E}[\omega]$, whatever $\tau$ and $\varepsilon$. Then, the worst case is always given by the function for $\overline{\tau}$ when $y<\mathbb{E}[\omega]$, and by the function for $\underline{\tau}$ when $y>\mathbb{E}[\omega]$. Finally, for cases~\ref{fig:DRnewsvendor-example1} and~\ref{fig:DRnewsvendor-example3}, the minimum value for the worst-case expected opportunity cost is located away from $y=\mathbb{E}[\omega]$, while for case~\ref{fig:DRnewsvendor-example2}, it is reached at  $y=\mathbb{E}[\omega]$.

As is expected when using distributionally robust optimization, we retrieve 2 interesting limit cases, i.e., for $\varepsilon=0$ and for a large-enough value of $\varepsilon$ (here, $\varepsilon = 1$ for the uniform ambiguity set, though potentially $\varepsilon\rightarrow \infty$ as $\rho \rightarrow 1$ for level-adjusted ambiguity sets). In that latter case, we end up with the solution to the robust optimization problem with uncertainty set $[0,1]$ for $\tau$.

\begin{corollary}
\label{corr:robnewsvendor2}
For sufficiently large values of $\varepsilon$, the radius of $\mathcal{B}_{\hat{\tau}} (\varepsilon)$, the solution of the distributionally robust Bernoulli newsvendor problem~\eqref{eq:DRnews2} converges to the robust solution $y^* \, = \, \mathbb{E}[\omega]$.
\end{corollary}

The proof of that corollary is given in Appendix~\ref{app:proofrobustnewsvendor2}. Eventually, our 2 limit cases are then the Bernoulli newsvendor solution $y^*=\hat{F}_\omega^{-1}(\hat{\tau})$ for $\varepsilon=0$, and the robust optimization outcome $y^*=\mathbb{E}[\omega]$ for sufficient large values of $\varepsilon$. In practice, the optimal value for $\varepsilon$ may be obtained in a data-driven framework, e.g., using cross-validation. \textcolor{black}{Alternatively, following the point of~\cite{Rahimian2019}, the ball radius in distributionally robust optimization may be chosen based on the risk profile of the market participant, thanks to the direct connection between distributionally robust optimization and optimization of coherent risk measures.}

\section{Simulation Study}
\label{sec:simul}

In order to illustrate the workings of the distributionally robust newsvendor problems described in the above, we perform a compact simulation study. Only the case of the ambiguity about $\hat{F}_s$ is covered, since the case of ambiguity about $\hat{F}_w$ is already discussed in the literature (especially in Refs~\cite{Fu2021} and~\cite{Rahimian2019}).

We illustrate some of the salient features of the distributionally robust approach to solving Bernoulli newsvendor problems based on a Monte-carlo simulation (with $N$ replicates), for a known distribution of the uncertain parameter $\omega$ and for a given $\tau$. For each replicate of the Monte-carlo experiment, we assume that it has been possible to observe $m$ previous outcomes of the Bernoulli process, and an estimate $\hat{\tau}$ of $\tau$ is obtained as an average of these $m$ estimates. Eventually, we will look at the expected loss $\mathbb{E}_\omega \left[ \mathcal{L}(y,\omega,\tau) \right]$ for the various approaches, i.e., direct solution to the Bernoulli newsvendor problem (Proposition~\ref{prop:Bnewsvendor}), the robust solution to the Bernoulli newsvendor problem (Corollary~\ref{corr:robnewsvendor2}) and the distributionally robust one (Theorem~\ref{theo:DRBnewsvendor2}) with both uniform and level-adjusted ambiguity sets. For the latter approach, an interesting aspect is to assess the impact of the ball radius on the outcome. An oracle is considered as benchmark, i.e., for the case where $\tau$ is known and not estimated. 

\begin{figure}[!ht]
    \centering
    \includegraphics[width=\columnwidth]{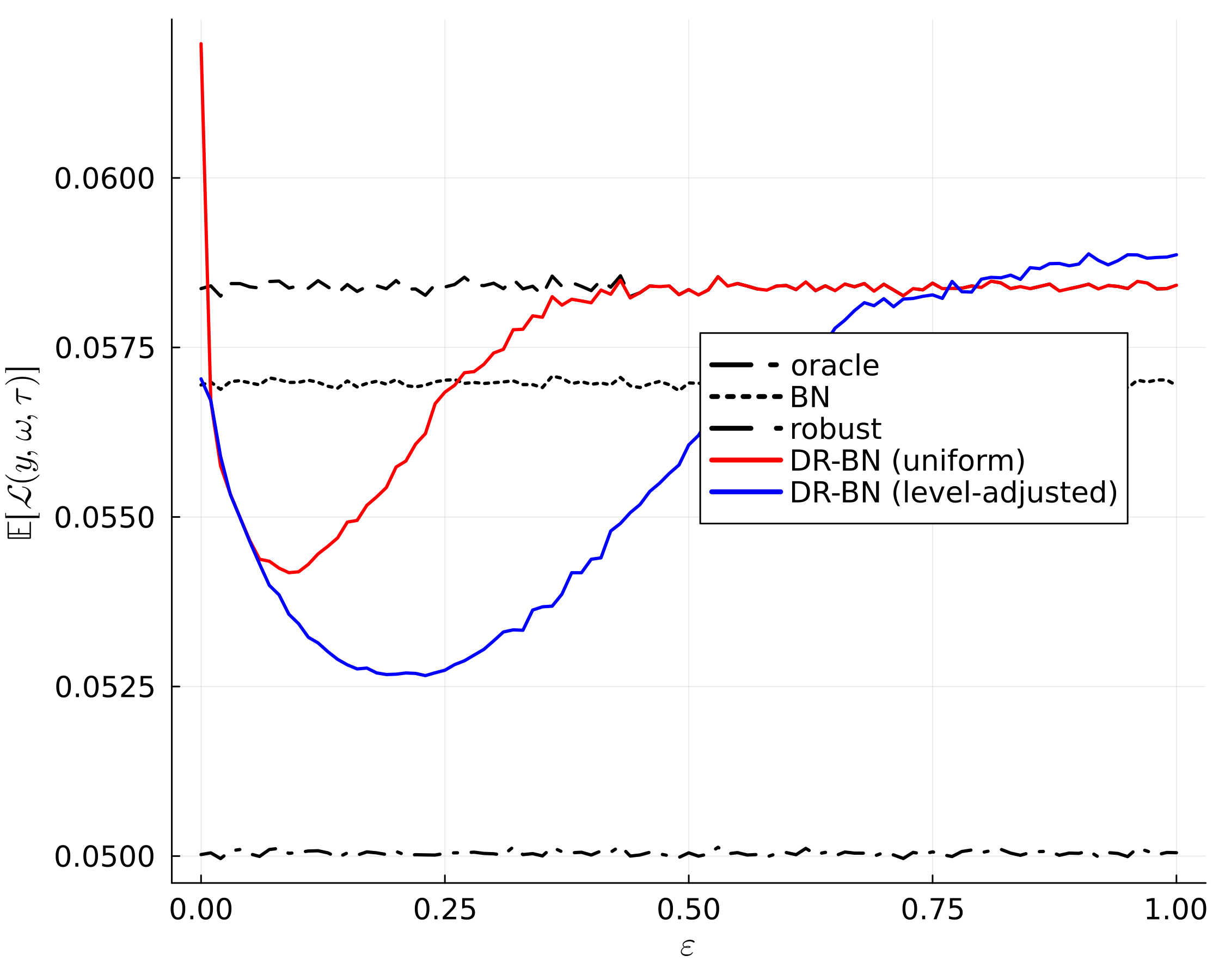}
    \caption{Expected loss $\mathbb{E}_\omega \left[ \mathcal{L}(y,\omega,\tau) \right]$ for the various approaches to solving the Bernoulli newsvendor problem (BN: Bernoulli Newsvendor; DR: Distributionally Robust). For the level-adjusted case, $\theta=0.9$.}
    \label{fig:DRvariedepsilon}
\end{figure}

For the sake of example, \textcolor{black}{we use $\omega \sim \text{Beta}(2,6)$ and $\tau=0.75$}. The results would be qualitatively similar if considering other distributions and values for $\tau$. The estimated chance of success $\hat{\tau}$ is obtained by taking the mean of $m=10$ random draws from a Bern($\tau$) distribution. \textcolor{black}{For that specific experiment,} expected loss function estimates are obtained based on $N=10^7$ replicates. Their evolution as a function of ball radius values, between 0 and 1, with an increment of $10^{-2}$, is depicted in Figure~\ref{fig:DRvariedepsilon}. 

As expected, the oracle yields the lowest expected loss since having perfect information about $\tau$. And, the expected loss for the direct and robust approach, as well as for the oracle, does not vary with $\epsilon$. It is the case for the distributionally robust approach, though, starting close to the direct solution for $\epsilon=0$ and eventually going towards the robust solution for larger values of $\epsilon$. There exists a range of values for $\epsilon$, for which the distributionally robust solution is the most competitive. In addition, this range is broader for the level-adjusted ambiguity set, while reaching a better minimum. \textcolor{black}{The case of $\varepsilon=0$ corresponds to the classical newsvendor problem solution, for which the uncertainty in the market outcomes (system length, as well as overage and underage penalties) is overlooked. The result of Figure~\ref{fig:DRvariedepsilon} hence clearly shows the benefits of accepting and accommodating the uncertainty in estimating these quantities, here all summarized within a single Bernoulli variable with chance of success $\tau$.} 

What the distributionally robust approach is aiming to do is to improve over the direct approach (using the estimated $\hat{\tau}$) and get closer to the oracle. Hence, a relevant performance measure can be defined as how much the distributionally robust approach improves over the direct approach, scaled by the difference between the direct approach and the oracle. We write $L_\text{BN}$ the expected loss for the direct approach, $L_\text{O}$ that for the oracle, and  $L^*_\text{DR-BN}$ the minimum one for the distributionally robust approach, given the optimal ball radius. The performance measure $\gamma$ is then defined as 
\begin{equation}
    \gamma = \frac{L_\text{BN} - L^*_\text{DR-BN}}{L_\text{BN} -  L_\text{O}} \, ,
\end{equation}
and can be expressed in percents. Let us denote by $\gamma_\text{U}$ and $\gamma_\text{LA}$ the performance measures for the uniform and level-adjusted ambiguity sets, respectively. For the example of Figure~\ref{fig:DRvariedepsilon}, we have $\gamma_\text{U} = 40.3$\% and  $\gamma_\text{LA} = 60.2$\%. When using real data, such a performance measure cannot be calculated, since the oracle is obviously not available. 

\begin{figure}[!ht]
    \centering
    \includegraphics[width=\columnwidth]{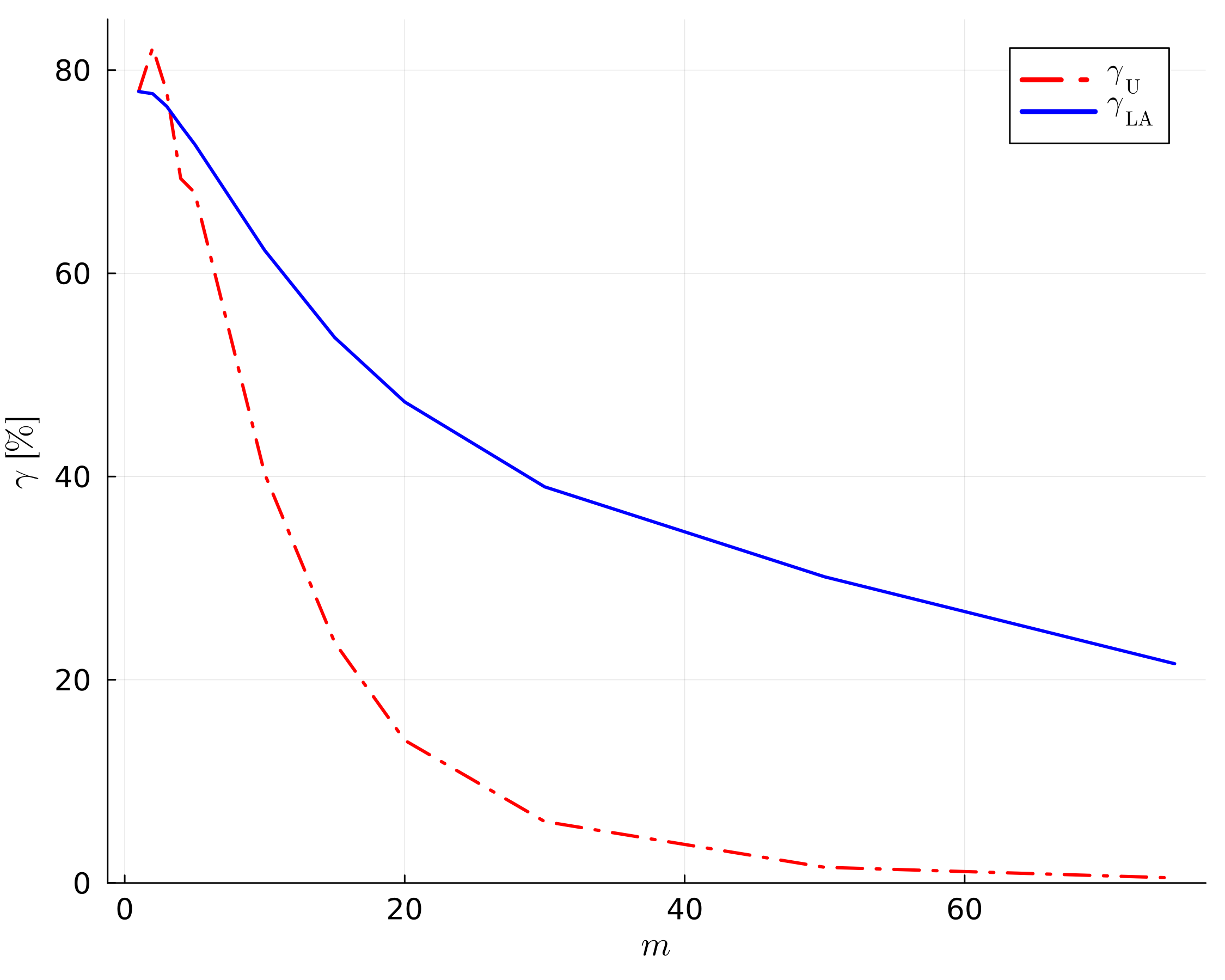}
    \caption{Expected value of $\gamma$ for both types of ambiguity sets (and $\theta=0.9$ for the level-adjusted one).}
    \label{fig:DRvariedm}
\end{figure}

The value of $\gamma$ is to vary with $\tau$, but most importantly it also varies with $m$ (being a proxy for the precision of the estimation of $\tau$). This is illustrated in Figure~\ref{fig:DRvariedm}, where $m$ ranges from 1 to 75 points used to estimate $\tau$, still for $\tau=0.75$. We also use $N =10^7$ replicates for that specific experiment.

As $m$ increases, the gap between the direct approach and the oracle reduces and the benefits from using the distributionally-robust approach with uniform ambiguity sets goes to 0 (this is why we stop at $m=75$). Eventually, if the estimation of $\tau$ was perfect, one would not need distributionally robust estimation anyway. The level-adjusted ambiguity sets allow to reach higher values for $\gamma$, whatever the precision in the estimation of $\tau$. For the uniform ambiguity sets, the benefits of the distributionally robust approach already drop close to 0 when reaching $m=40$. Here, a value of $\theta=0.9$ was used for the sake of illustration. For a real-world case study, one would need to decide upon $\theta$ through cross-validation, in addition to the ball radius $\varepsilon$.

\section{Case-study Application}
\label{sec:appl}

\subsection{Dataset and Experimental Setup}

Emphasis is placed on a case-study with real data in France. Wind power generation data is available for a portfolio of wind farms in the mid-West of France, the location of which cannot be disclosed for confidentiality reasons. The power generation from that portfolio is normalized by its nominal capacity, so that all values are within the unit interval. Data is available for a period of 2 years, from October 2018 to September 2020 (both months included), i.e., exactly 731 days. We therefore look at the case of this portfolio of wind farms participating in the French electricity market, with day-ahead market data obtained from EEX and balancing market data obtained from the French Transmission System Operator (RTE). \textcolor{black}{Since the French electricity market switched from a two-price to a one-price imbalance settlement during that period, while the approach proposed here is for the two-price imbalance settlement case, the data is adapted accordingly: in terms of overage and underage penalties, this means that all cases where such penalties could be negative are removed and rounded to 0 instead.} It is not a very penalizing market, with the spread between forward and balancing prices being on average of 13.5\% of the forward prices.

For wind power generation, probabilistic forecasts are generated for that same period. Forecasts are issued every day before gate closure for the following day. They are generated with a pattern matching approach, in a fashion similar to that described in \cite{vanderMeer2022}\footnote{Another approach to generate probabilistic forecasts was also employed, based on quantile regression forests and directly using the functions in \cite{Meinshausen2017}. It led to similar qualitative results from trading, even though forecasts and their quality were different.}. The scenarios produced with that method originally are used to deduce predictive cumulative distributions functions, summarized by quantiles with levels from 0.025 to 0.975, with 0.05 increments. The quality of the forecasts was assessed and deemed in par with the state of the art. Especially, these forecasts are probabilistically calibrated, hence with the observed coverage of the quantiles being consistent with the nominal levels expected. It is an important aspect to insure that no systematic error leads to biases in decision-making.

In terms of the chance of success of the Bernoulli variables, a separate model is used for each hour of the day, in order to capture daily patterns that are commonly present in electricity markets. However, the models employed are kept simple: for a given hour of the day, a forecast is produced as a moving average, i.e., by taking the average of overage and underage penalties over the last $m$ days. $m$ then is a parameter than can be optimized within a cross-validation framework. More advanced models were benchmarked, e.g., based on exponential smoothing and Holt-Winters generalizations \cite{Jonsson2014}, though they did not lead to improved outcomes.

\textcolor{black}{As is usually the case for data-driven approaches, the available data is split into a number of subsets. Here, the first 131 days are used as a warm start (since needing an initial window to obtain an estimate for $\tau$ based on observed market outcomes -- here, 91 days) and for cross-validation (here, 40 days).} Then, the remaining 600 days are used for genuine out-of-sample evaluation of the trading strategies and observed outcomes. Cross-validation is used in order to decide upon values for the ball radius $\rho$ (for the case of ambiguity about $\hat{F}_\omega$) and $\varepsilon$ (for the case of ambiguity about $\hat{F}_s$), as well as the shape parameter $\theta$ (for the level-adjusted ball). In addition, the value of $m=90$ days was found to yield the best results.

\textcolor{black}{In practice, we considered two different approaches consisting of \emph{(i)} using a single and fixed period for cross-validation, \emph{(ii)} using a sliding window for cross-validation, allowing for temporal variations in the optimal parameters. In both cases, the length of the cross-validation period is of 40 days. Comparable results were obtained in both cases, hence only the results for selection based on a single and fixed period are presented here.} We end up with 3 setups, since for the case of ambiguity about $\hat{F}_s$, we use both uniform and level-adjusted ambiguity sets. Final parameters are gathered in Table~\ref{tab:parameters}.

\begin{table}[!ht]
    \centering
    \caption{Parameters for the distributionally robust approaches. \label{tab:parameters}}
    \begin{tabular}{l c c c}
    \hline
    \hline
     Ambiguity on & $\hat{F}_\omega$ & $\hat{F}_s$ & $\hat{F}_s$\\
     Ball radius & $\rho=0.24$  & $\varepsilon=0.15$ &  $\varepsilon=0.25$ \\
     Shape parameter & - & - & $\theta=0.5$ \\
    \hline
    \hline
    \end{tabular}
\end{table}

\subsection{Results and Discussion}

Two approaches are considered as benchmarks, i.e., the oracle consisting in having perfect forecast for renewable energy generation, and the direct approach to solving Bernoulli newsvendor problems (as in Proposition~\ref{prop:Bnewsvendor}). For all approaches, we look at the revenue $R$ obtained per MWh generated as well as the regret $r$ per MWh generated. The regret is defined as the difference between the revenue of the oracle $R_{\text{O}}$ and that of the approach considered, e.g., $R_{\text{BN}}$ for the direct newsvendor approach, $r_{\text{BN}} = R_{\text{O}} - R_{\text{BN}}$. Revenues and regret values for the various approaches are collated in Table~\ref{tab:results}. We additionally give there the advantage ratio for the various methods, i.e., the percentage of times the revenue they yield is at least as good as that obtained with the direct approach.

\begin{table}[!htbp]
    \centering
    \caption{Revenues, regret and advantage ratio for the various approaches, over the 600-day evaluation period. \label{tab:results}}
    \begin{tabular}{l c c c}
    \hline
    \hline
    Approach & $R$ [\euro/MWh]& $r$ [\euro/MWh] & Adv.~ratio [\%]  \\
    \hline
    Oracle & 31.63 & 0 & 100 \\
    BN & 29.25 & 2.38 & - \\
    DR-BN (amb. $\hat{F}_w$) & 29.34 & 2.29 & 76.52 \\
    DR-BN & 29.35 & 2.28 & 74.16 \\
    (amb. $\hat{F}_s$, uniform) & & & \\
    DR-BN & 29.36 & 2.27 & 74.46\\
    (amb. $\hat{F}_s$, level-adjusted) & & & \\
    \hline
    \hline
    \end{tabular}
\end{table}

\begin{figure}[!ht]
    \centering
    \includegraphics[width=\columnwidth]{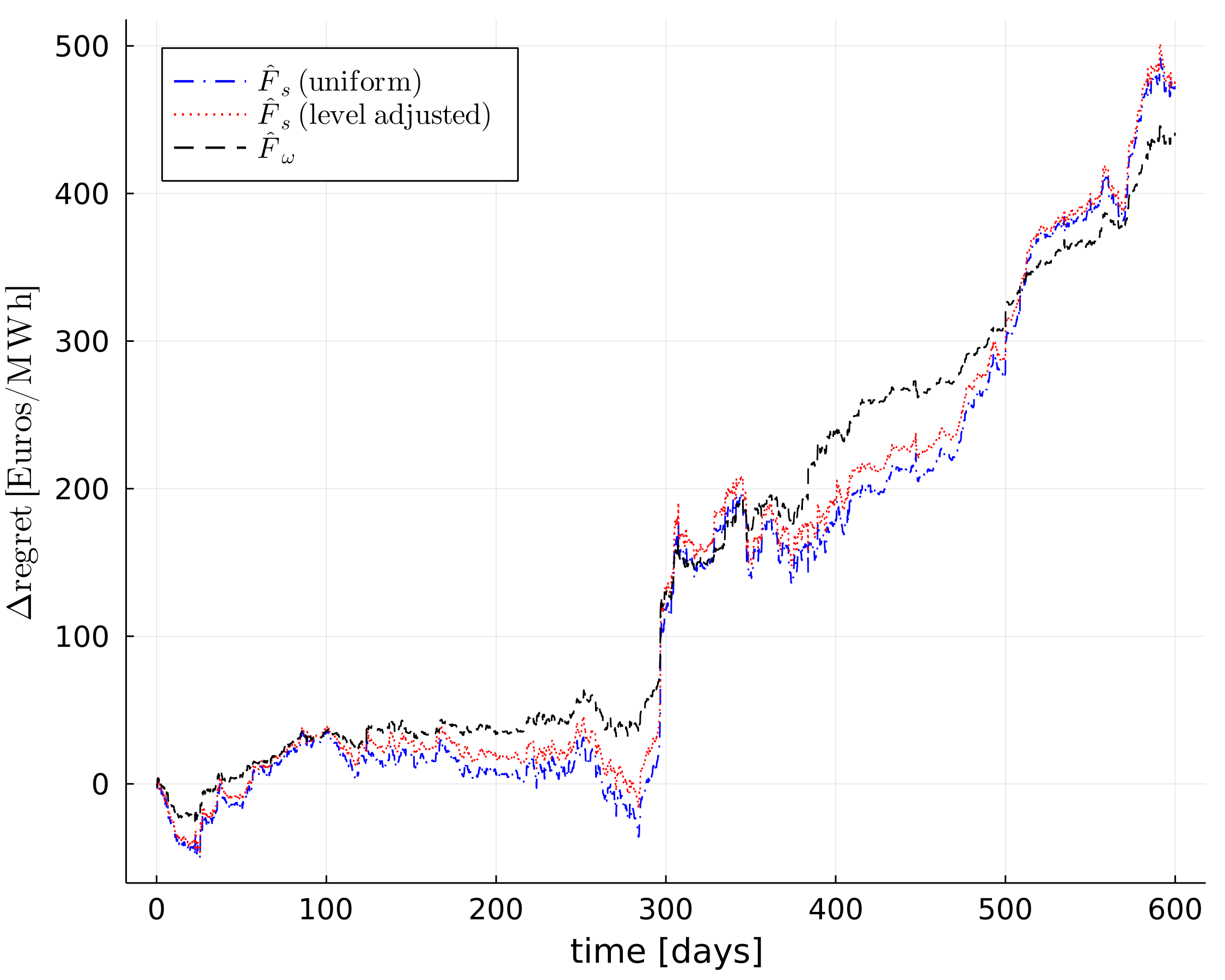}
    \caption{Cumulative difference in regret (with respect to direct approach to solving the Bernoulli newsvendor problem) for the various distributionally robust optimization approaches.}
    \label{fig:cumdregret}
\end{figure}

The improvements for the various distributionally robust opimization approaches are modest, but still ranging from 3.8\% (for the ambiguity about $\hat{F}_w$) to 4.6\% (for the ambiguity about $\hat{F}_s$ and the level-adjusted ambiguity set) in terms of lowering the regret. And, these improvements come for free, nearly, since only resulting from updated formulas to solve the Bernoulli newsvendor problem. In addition, these improvements are fairly consistent over time. This is firstly supported by values for the advantage ratio in Table~\ref{tab:results}, with a slight advantage for the distributionally approach with ambiguity about $\hat{F}_\omega$. It is also supported by Figure~\ref{fig:cumdregret}, which shows the cumulative difference in regret (denoted $\Delta \text{regret}$) between the direct approach and the various distributionally robust approaches. There, positive values mean that the distributionally robust approach does better than the direct approach, while an increasing trend tells that the former does consistently (in time) better than the latter.

In this figure, it is interesting to see that there are different periods over which the approaches considering ambiguity about $\hat{F}_w$ and $\hat{F}_s$ outperform each other. \textcolor{black}{This shows} that accommodating the lack of forecast quality for probabilistic forecasts of wind power generation and market penalties may have more or less relevance over successive periods. Most importantly, it hints at the fact it may make sense to generalize the approach to jointly consider ambiguity about $\hat{F}_w$ and $\hat{F}_s$ in the future. Finally, even if these differences may appear modest, they still translate to non-negligible monetary differences. Taking the example of a 100 MW wind farm with a capacity factor of 40\% for instance, earning 11 c\euro\ more per MWh produced yields an additional revenue of more nearly 40 000 \euro\ per year. For the case of this French dataset, the market penalties are low on average. However, with larger penalties (for instance experimented here by multiplying market penalties by a constant factor), the benefits from employing distributionally robust approaches would increase accordingly.

\section{Conclusions}
\label{sec:concl}

After reformulating the problem of renewable energy offering in electricity markets as a Bernoulli newsvendor problems, various alterations were considered. Especially, emphasis is placed on the fact that forecasts for renewable energy generation and market quantities are not perfect, and this can be accommodated using distributionally robust optimization. While the solution of distributionally robust Bernoulli newsvendor problem can be deduced from some recent results in the literature, a closed-form solution is provided here for the case of ambiguity about the Bernoulli random variable.

Future works should focus on jointly accommodating ambiguity about the forecasts for renewable energy generation and the chance of success of the Bernoulli random variable. And, since such approaches are data-driven and with parameters to be decided upon, it may be relevant to look at automatic ways to optimize them in an online environment. In addition, the results should be reproduced in different markets and with various types of forecasts. Finally, the distributionally robust approach should be adapted to the case of one-price imbalance settlement, which was not considered here.

\section*{Acknowledgments}

Acknowledgements are due to Dennis van der Meer and Simon Camal at Mines Paris (France) for producing (and sharing) the wind power forecasts data, Vestas for providing wind power data, ECMWF the input weather forecasts, and EEX the power market data. Finally, Gustav Klenz Larsen and Jalal Kazempour at DTU are acknowledged for many relevant discussions about distributionally robust offering strategies for renewable energy producers. 4 reviewers and an associate editor are also to be acknowledged for their comments and suggestions, which allowed to greatly improved the version of the manuscript submitted originally.

\appendix

\subsection{Proof of Proposition~\ref{prop:Bnewsvendor}}
\label{app:proofnewsvendor}

For the more general case of the Bernoulli newsvendor problem, the expected opportunity cost relies on 2 random variables $\omega$ and $s$. It is defined as
\begin{align}
    \mathbb{E}_{\omega,s} \left[ \mathcal{L}(y,\omega,s) \right] = &  \, \, \int_s \int_0^y (y-\omega) f_\omega(\omega) f_s(s) d\omega ds \nonumber \\ & + \,  \int_s \int_y^1 (\omega-y)f_\omega(\omega) f_s(s) d\omega ds \, . \nonumber \nonumber
\end{align}
Based on the price-taker assumption, the random variables $\omega$ and $s$ are independent, allowing us to obtain
\begin{align}
    \mathbb{E}_{\omega,s} \left[ \mathcal{L}(y,\omega,s) \right] = & \, \,  \text{P}[s=0] \int_0^y (y-\omega) f_\omega(\omega) d\omega  \nonumber \\ & + \,  \text{P}[s=1] \int_y^1 (\omega-y)f_\omega(\omega) d\omega  \, . \nonumber 
\end{align}
This readily gives
\begin{align}
    \mathbb{E}_{\omega,s} \left[ \mathcal{L}(y,\omega,s) \right] = & \, \,  (1-\tau) \int_0^y (y-\omega) f_\omega(\omega) d\omega  \nonumber \\ & + \, \tau \int_y^1 (\omega-y)f_\omega(\omega) d\omega \nonumber
\end{align}
One recognizes that this is the same problem as for the classical newsvendor problem. To be self consistent, we describe the solution approach in the following.

The expected opportunity cost is a convex function defined on a compact set. One can find its minimum by looking for the zero of its derivative. By differentiating using the Leibniz rule (differentiation under the integral sign), one has
\begin{subequations}
\begin{align}
    \frac{\partial }{\partial y} \int_0^y (1-\tau) (y-\omega)f_\omega(\omega) d\omega  &= \int_0^y \frac{\partial }{\partial y} (1-\tau) (y-\omega)f_\omega(\omega) d\omega \nonumber \\ 
    & = (1-\tau) \int_0^y f_\omega(\omega) d\omega \ , \nonumber
\end{align}
and
\begin{align}
\frac{\partial }{\partial y} \int_y^1 \tau (\omega-y)f_\omega(\omega) d\omega & = \int_y^1 \frac{\partial }{\partial y} \tau (\omega-y)f_\omega(\omega) d\omega \nonumber \\
& = - \tau \int_y^1 f_\omega(\omega) d\omega \ . \nonumber
\end{align}
\end{subequations}
Consequently,
\begin{subequations}
\begin{align}
  \frac{\partial }{\partial y} \, \mathbb{E}_\omega \left[ \mathcal{L}(y,\omega;\tau) \right] & = (1-\tau) \int_0^y f_\omega(\omega) d\omega - \tau \int_y^1 f_\omega(\omega) d\omega \nonumber \\
  & = (1-\tau)  F_\omega (y) - \tau [1 - F_\omega (y) ] \nonumber
\end{align}
\end{subequations}
Finally,
\begin{equation}
    \frac{\partial }{\partial y} \, \mathbb{E}_\omega \left[ \mathcal{L}(y^*,\omega;\tau) \right]  = 0 \quad \implies \quad y^* = F_\omega^{-1} \left( \tau \right) \nonumber
\end{equation}

\subsection{Proof of Corollary~\ref{corr:robnewsvendor1}}
\label{app:proofrobustnewsvendor1}

Theorem~\ref{theo:DRBnewsvendor1} tells us that the solution of the distributionally robust newsvendor problem (with ambiguity about $\hat{F}_\omega$) is given by~\eqref{eq:DROsol1}. The robust limiting case is obtained for $\rho=1$. Based on the definition of deformation operators in Definition~\ref{def:deformationoperators}, we know that in such a case 
\begin{equation}
    \overline{F}_\omega = H(0), \quad \underline{F}_\omega = H(1). \nonumber
\end{equation}
Consequently, whatever the value of $\hat{\tau}$ we can deduce that
\begin{equation}
    \overline{F}_\omega{}^{-1}(\hat{\tau}) = 0, \quad \underline{F}_\omega{}^{-1}(\hat{\tau}) = 1. \nonumber
\end{equation}
Plugging these values in~\eqref{eq:DROsol1} yields $y^*=\hat{\tau}$.

\subsection{Proof of Theorem~\ref{theo:DRBnewsvendor2}}
\label{app:proofDRnewsvendor2}

The proof is a bit more involved than for the classical newsvendor problem, since we first need to reformulate the maximisation problem (over $\tau \in \mathcal{B}_{\hat{\tau}}(\varepsilon)$) and then find the solution of the resulting minimization problem. 

\subsubsection{Reformulation of the Maximisation Problem}

In this first part, we are looking into reformulating the inner maximisation problem, over $\tau \in \mathcal{B}_{\hat{\tau}} (\varepsilon)$. We first assess whether $\mathbb{E}_{\omega,s} [\mathcal{L}(y,\omega,s)]$ is monotonous in $\tau$. We have, for a given $\delta>0$,
\begin{align}
\mathbb{E}_{\omega,s}[\mathcal{L}(y,\omega,\tau+\delta)] \, &  = \, \int_0^y (1-(\tau+\delta)) (y-\omega) f_\omega(\omega) d\omega \nonumber \\ & + \, \int_y^1 \tau+\delta (\omega-y) f_\omega(\omega) d\omega \nonumber \\
& = \, \mathbb{E}_{\omega,s}[\mathcal{L}(y,\omega,\tau)] - \int_0^y \delta (y-\omega) f_\omega(\omega) d\omega \nonumber \\ & +  \, \int_y^1 \delta (\omega-y) f_\omega(\omega) d\omega \nonumber \\
& = \, \mathbb{E}_{\omega,s}[\mathcal{L}(y,\omega,\tau)] + \delta \int_0^y (\omega-y) f_\omega(\omega) d\omega \nonumber \\ & + \, \int_y^1 \delta (\omega-y) f_\omega(\omega) d\omega \nonumber \\
& = \, \mathbb{E}_{\omega,s}[\mathcal{L}(y,\omega,\tau)] + \delta \int_0^1 \omega f_\omega(\omega) d\omega \nonumber \\ & - \, \delta \, y \, \int_0^1 f_\omega(\omega) d\omega \nonumber
\end{align}
We therefore obtain that 
\begin{equation}
\mathbb{E}_{\omega,s}[\mathcal{L}(y,\omega,\tau+\delta)] - \mathbb{E}_{\omega,s}[\mathcal{L}(y,\omega,\tau)] = \delta \left( \mathbb{E}[\omega] - y \right) \nonumber
\end{equation}
which shows that $\mathbb{E}_{\omega,s}[\mathcal{L}(y,\omega,\tau)]$ is monotonous in $\tau$. It is strictly (and linearly) increasing for $y\leq\mathbb{E}[\omega]$, and then strictly (and linearly) decreasing for $y>\mathbb{E}[\omega]$.

By extension, we can deduce that the worst cases are at the boundary of the ball $\mathcal{B}_{\hat{\tau}} (\varepsilon)$. For $y\leq\mathbb{E}[\omega]$, the worst case is for $\overline{\tau} = \hat{\tau} +\varepsilon$, and for $y>\mathbb{E}[\omega]$, the worst case is for $\underline{\tau} = \hat{\tau} - \varepsilon$. We can then reformulate the inner maximization problem as
\begin{align}
    \max_{\tau \in \mathcal{B}_{\hat{\tau}}(\varepsilon)} \, \mathbb{E}_{\omega,s}[\mathcal{L}(y,\omega,\tau)] \, \equiv & \, \,  \mathbb{E}_{\omega}[\mathcal{L}(y,\omega,\overline{\tau})] \mathbf{1}_{\{y\leq\mathbb{E}[\omega]\}} \nonumber \\ & + \mathbb{E}_{\omega}[\mathcal{L}(y,\omega,\underline{\tau})] \mathbf{1}_{\{y>\mathbb{E}[\omega]\}} \, . \nonumber
\end{align}

\subsubsection{Solution of the Minimisation Problem}

Based on the above, the min-max problem can be rewritten as
\begin{equation}
    y^* \, = \,  \argmin_y \, \, \mathbb{E}_{\omega}[\mathcal{L}(y,\omega,\overline{\tau}] \,  \mathbf{1}_{\{y\leq\mathbb{E}[\omega]\}} + \mathbb{E}_{\omega}[\mathcal{L}(y,\omega,\underline{\tau})] \,  \mathbf{1}_{\{y>\mathbb{E}[\omega]\}} \, . \nonumber
\end{equation}

Since we necessarily have $\hat{F}_\omega^{-1}(\overline{\tau}) \geq \hat{F}_\omega^{-1}(\underline{\tau}) $, it is not possible to have both $\hat{F}_\omega^{-1}(\overline{\tau})  \leq \mathbb{E}[\omega]$ and $\hat{F}_\omega^{-1}(\underline{\tau}) > \mathbb{E}[\omega] $. Therefore, the solution of the above problem is given by
\begin{align}
    y^* \, = & \, \,   \hat{F}_\omega^{-1}(\overline{\tau}) \,   \mathbf{1}_{\{ \hat{F}_\omega^{-1}(\overline{\tau}) <\mathbb{E}[\omega]\}} + \hat{F}_\omega^{-1}(\underline{\tau}) \,  \mathbf{1}_{\{ \hat{F}_\omega^{-1}(\underline{\tau})  >\mathbb{E}[\omega]\}} \nonumber \\
     & + \, \mathbb{E}[\omega] \, \mathbf{1}_{\{ \hat{F}_\omega^{-1}(\overline{\tau}) \geq \mathbb{E}[\omega]\}} \, \mathbf{1}_{\{ \hat{F}_\omega^{-1}(\underline{\tau}) \leq \mathbb{E}[\omega]\}} \, . \nonumber
\end{align}

\subsection{Proof of Corollary~\ref{corr:robnewsvendor2}}
\label{app:proofrobustnewsvendor2}

In the robust optimization case (with ambiguity about $\hat{F}_s$), the min-max optimization problem is
\begin{equation}
    y^* = \argmin_{y} \, \max_{\tau \in [0,1] } \, \mathbb{E}_{\omega,s} \left[ \mathcal{L}(y,\omega,\tau) \right] \, . \nonumber
\end{equation}
As for the distributionlly robust case, we end up with a worst-case function which is a piece-wise function
\begin{align}
    \max_{\tau \in [0,1]} \, \mathbb{E}_{\omega,s}[\mathcal{L}(y,\omega,\tau)] \, \, \equiv & \, \, \mathbb{E}_{\omega}[\mathcal{L}(y,\omega,1)] \mathbf{1}_{\{y\leq\mathbb{E}[\omega]\}} \nonumber \\ & + \mathbb{E}_{\omega}[\mathcal{L}(y,\omega,0)] \mathbf{1}_{\{y>\mathbb{E}[\omega]\}} \, . \nonumber
\end{align}
Since that piece-wise function is decreasing for $y<\mathbb{E}[\omega]$ and then increasing for $y>\mathbb{E}[\omega]$, its minimum is at $y^* = \mathbb{E}[\omega]$.

\vfill
\begin{IEEEbiography}[{\includegraphics[width=1in,height=1.25in,clip,keepaspectratio]{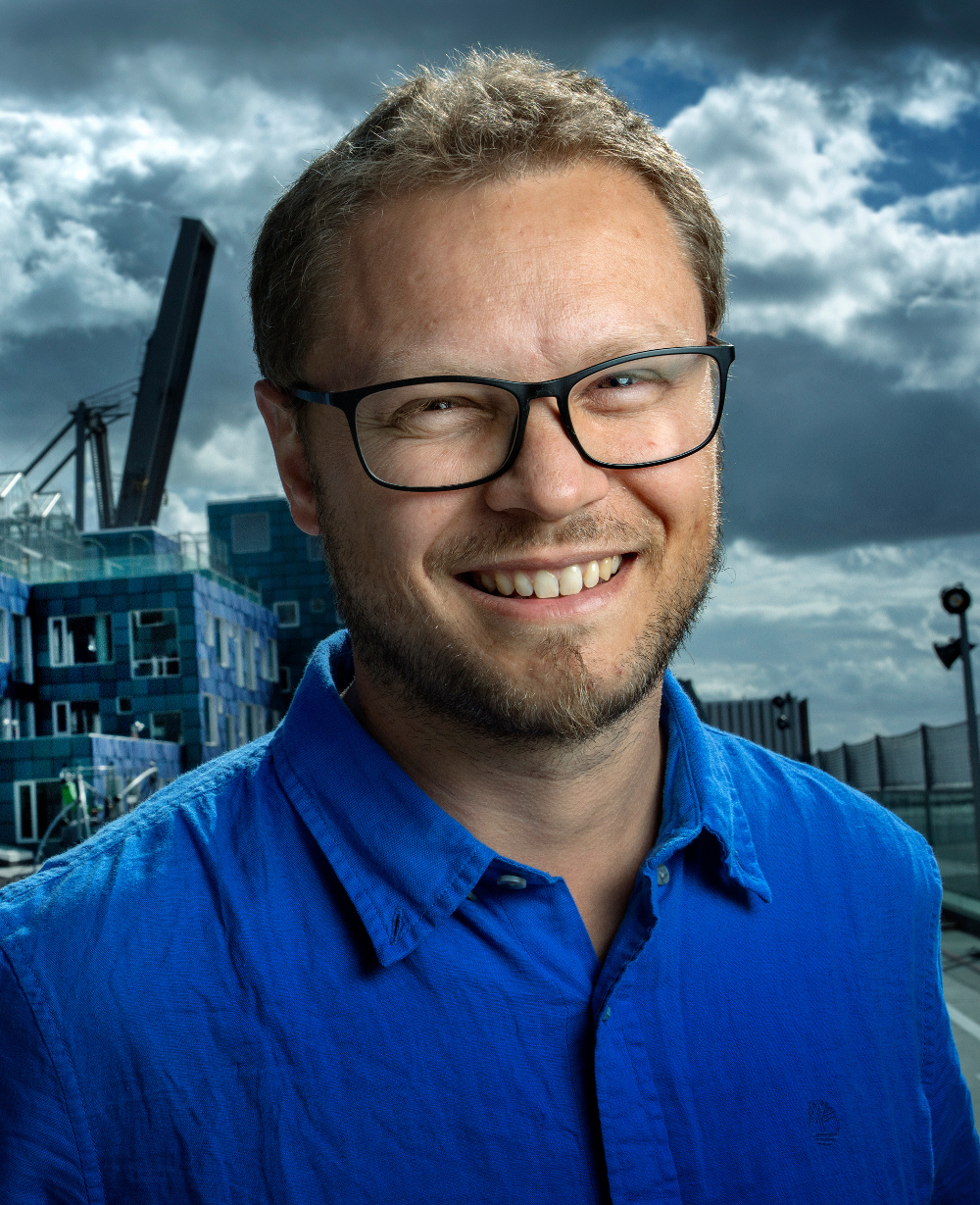}}]{Pierre Pinson} (Fellow, IEEE) received the M.Sc. degree in applied mathematics from the National Institute of Applied Sciences (INSA), Toulouse, France, in 2002 and the Ph.D. degree in energetics from Ecole des Mines de Paris, France, in 2006. He is the chair of data-centric design engineering at Imperial College London, United Kingdom, Dyson School of Design Engineering. He is also an affiliated professor of operations research and analytics with the Technical University of Denmark and a chief scientist at Halfspace (Denmark). He is the editor-in-chief of the \emph{International Journal of Forecasting}. His research interests include analytics, forecasting, optimization and game theory, with application to energy systems mostly, but also logistics, weather-driven industries and business analytics. He is a Fellow of the IEEE, an INFORMS member and a director of the International Institute of Forecasters (IIF).
\end{IEEEbiography}

\end{document}